# Old and new geometric polyhedra with few vertices

S. Lawrencenko and A. Lao

**ABSTRACT**

This paper deals with triangulations of the 2-torus with the vertex labeled general octahedral graph $O_4$ which is isomorphic to the complete four-partite graph $K_{2,2,2,2}$; it is known that there exist precisely twelve such triangulations. We find all the 12 triangulations in a Schlegel diagram of the hyperoctahedron and realize all of them geometrically with the same 1-skeleton in 3-space. In particular, we identify two geometric polyhedral tori (both without self-intersections) with the same 1-skeleton in 3-space, but without a single common face, or in other words their intersection (as point-sets) is only their common 1-skeleton. Similarly, all the twelve triangulations of the 2D projective plane with the vertex labeled complete graph $K_6$ are found in a Schlegel diagram of the 5-simplex and all are realized geometrically with the same 1-skeleton in 4-space; especially we obtain a pair of triangulations of the Möbius band and a pair of triangulated projective planes with the same 1-skeleton (within each pair) in 3-space and 4-space, respectively, without a single common face. The constructed polyhedra are modeled and visualized with GeoGebra.

## 1 INTRODUCTION

A *geometric realization* of an abstract triangulation $T$ of a closed $d$-manifold is a geometric polyhedron (with all faces represented by geometric simplices) in $n$-space $\mathbf{R}^n$ ($d \leq n-1$), whose adjacencies correspond to the adjacencies of $T$ or in other words whose polyhedral structure is *isomorphic* to $T$ as triangulations.

It should be emphasized that the problem of finding a geometric realization of a given abstract triangulation is not always an easy one since the realizing polyhedron is required to be free from self-intersections or in other words the realizing polyhedron should provide a geometric *embedding* of $T$ in $\mathbf{R}^n$. By now some results are available for $d=2$; in particular, it is known that every triangulation of the (2-)sphere [18] or (2-)torus [1] is geometrically realizable in $\mathbf{R}^3$ while every triangulation $T$ of the (2D) projective plane [2] is geometrically realizable in $\mathbf{R}^4$ and moreover [3] has a face $f$ such that the triangulation of the Mobius band obtained from $T$ by removing the interior of $f$ has a geometric realization in $\mathbf{R}^3$. In this paper we focus on geometric realization of two specific abstract triangulations. One is the 6-regular triangulation of the torus with 8 vertices shown in Figure 1 (left), the other is the 5-regular triangulation of the projective plane with the complete graph $K_6$ shown in Figure 1 (right); identify the sides of each fundamental polygon in pairs, as prescribed by the labels, to obtain a torus or projective plane, respectively.

The triangulation shown in Figure 1 (left) is realized geometrically [14] as a toroidal polyhedral suspension in $\mathbf{R}^3$ and also as a 2D *noble* polyhedron in 4-space $\mathbf{R}^4$, that is, a polyhedron which is *isohedral* (that is, all faces are similar) and *isogonal* (that is, all vertices are similar), whose properties are studied extensively in [17].

The importance and significance of the triangulation of the torus shown in Figure 1 (left) is justified by the fact that its graph is isomorphic to the complete four-partite graph $K_{2,2,2,2}$ obtained by removing a perfect matching from the complete graph $K_8$. The graph $K_{2,2,2,2}$ is also known as a special case of the *general octahedral graph* $O_d$, for dimension $d = 4$, since it is the graph of the *hyperoctahedron*, a.k.a. 4D regular *cross-polytope* (in $\mathbf{R}^4$) or the *16-cell* or 4D *orthoplex*. The hyperoctahedron is a convex regular 4-polytope bounded by sixteen 3D facets (that is, 3-faces which are geometric regular 3-simplices, or tetrahedra). The hyperoctahedron has 8 vertices, 24 edges, and 32 triangular (2-)faces. Geometrically, the hyperoctahedron can be imagined as the convex hull of some four pairwise perpendicular intervals in $\mathbf{R}^4$, intersecting at their common midpoint. The hyperoctahedron is one of the six regular convex 4-polytopes in $\mathbf{R}^4$ (a.k.a. *polychora*). The hyperoctahedron is a natural generalization of the 3D octahedron (bounded by eight 2-simplexes) which in its turn is a generalization of the 2D rhombus (bounded by four 1-simplexes); in this context, the usual octahedron is understood as an 8-cell and the rhombus as a 4-cell.

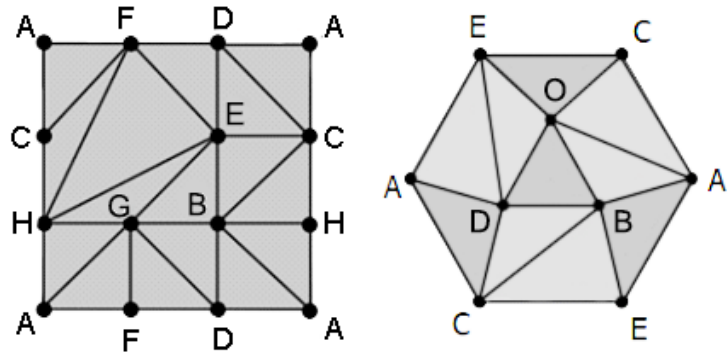

Figure 1: Abstract triangulations of the torus (left) and projective plane (right)

The triangulation in Figure 1 (left), respectively Figure 1 (right), is found in the 2-skeleton of the hyperoctahedron (a regular convex 4-polytope in $\mathbf{R}^4$), respectively the 5-simplex (a regular convex 5-polytope in $\mathbf{R}^5$). This observation leads us (Sections 4, 5) to new polyhedra found in the Schlegel diagrams of the convex polytopes in $\mathbf{R}^3$ and $\mathbf{R}^4$, respectively. A *Schlegel diagram* is obtained by a central projection of the polytope onto one of its facets.

For surfaces other than the sphere, the following *flexibility phenomenon* may be the case: There may exist different triangulations of the same surface with the same vertex labeled graph, where two triangulations are considered different provided that they have different face sets. Especially it is shown that both the torus [8, 10, 15] and the projective plane [5, 6, 9, 12] admit precisely 12 pairwise different triangulations with the vertex labeled graphs $O_4 = K_{2,2,2,2}$ and $K_6$ respectively (which is a coincidence). For example, the triangulation in Figure 2 (left) is identical with the initial one in Figure 1 (left) but different from the one in Figure 2 (right). The flexibility phenomenon is never the case in the sphere since by a theorem of Whitney [19] each 3-connected planar graph G has a unique dual and thus is uniquely embeddable in the sphere.

Furthermore [7, 15], the 12 triangulations of each named surface split into 6 pairs of mutually complementary triangulations. Two triangulations $T_1$ and $T_2$ of the same surface with the same vertex labeled graph $G$ are called *complementary* to each other if each cycle of three edges of $G$ bounds a face in $T_1$ if and only if it does not bound a face in $T_2$. The intersection of the face sets of mutually complementary triangulations $T_1$ and $T_2$ is empty; check with Figure 2 which presents a pair of mutually complementary triangulations of the torus with the vertex labeled graph $O_4 = K_{2,2,2,2}$.

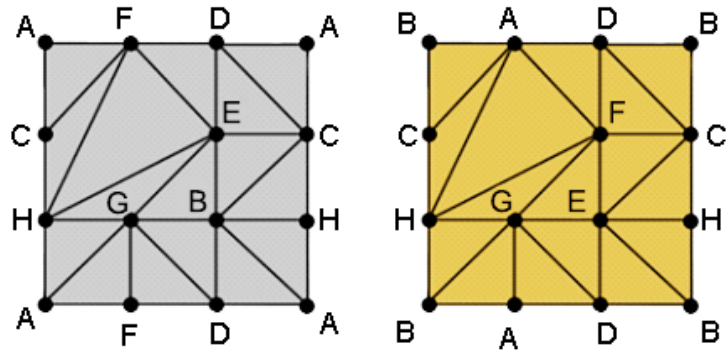

Figure 2: The initial triangulation of the torus (left) and its complementary triangulation (right)

It should be stressed that we realize geometrically *all* the 12 triangulations of the torus, respectively projective plane, using as the 1-skeleton the same geometric realization of the graph $O_4 = K_{2,2,2,2}$, respectively $K_6$, found in $\mathbf{R}^3$, respectively $\mathbf{R}^4$. In the torus case (Section 4), the construction provides an example of two polyhedral tori (both without self-intersections) having a common 1-skeleton in 3-space but having no face in common; this example justifies the epigraph of the current paper.

## 2 THE TOROIDAL POLYHEDRAL BIPYRAMIDAL SUSPENSION

In this section we describe the first author's earlier construction [14] which geometrically realizes the torus triangulation in Figure 1 (left) as a toroidal polyhedral *(bipyramidal) suspension*, that is, with the property that all but two vertices of $O_4 = K_{2,2,2,2}$ lie in the same plane, $z = 0$.

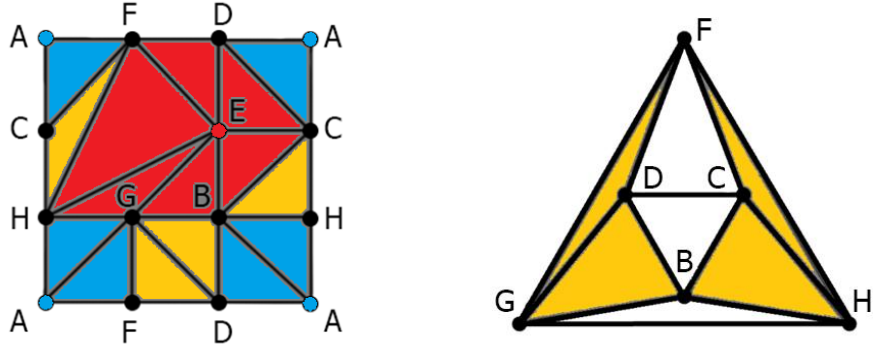

Figure 3: The torus triangulation (left) and its equatorial plane 2-complex (right)

For convenience of observation, the diagram in Figure 1 (left) is redrawn in Figure 3 (left) with colors. The crucial observation is that by judiciously removing only two vertices $A$ ("north pole") and $E$ ("south pole") from the diagram in Figure 3 (left) together with the (blue) faces incident with $A$ and the (red) faces incident with $E$, we obtain a planar graph and a planar simplicial 2-complex in the $xy$-plane thought of as the "equatorial plane"; see Figure 3 (right). To construct a toroidal polyhedral suspension, we assign coordinates to the vertices of the abstract triangulation in Figure 3 (left) as follows:

$$A = (0,0,\sqrt{8}), \quad B = \left(-\frac{\sqrt{8}}{k}, 0, 0\right), \quad C = \left(\frac{\sqrt{2}}{k}, \frac{\sqrt{6}}{k}, 0\right), \quad D = \left(\frac{\sqrt{2}}{k}, -\frac{\sqrt{6}}{k}, 0\right),$$
$$E = (0,0,-\sqrt{8}), \quad F = (\sqrt{8},0,0), \quad G = (-\sqrt{2},-\sqrt{6},0), \quad H = (-\sqrt{2},\sqrt{6},0),$$

where $k$ is a parameter which can be varied in the model ($2 < k < \infty$). We regard this parameter as the absolute reciprocal value of the negative coefficient of homothety, the center of homothety is the common circumcenter of the nonfacial regular triangular cycles $BCD$ and $FGH$; the homothety maps $FGH$ onto $BCD$ (check with Figure 3, right). The poles $A$ and $E$ are placed in the $z$-axis above and below the equatorial plane, respectively. The quadruples of vertices $A, F, G, H$ and $E, F, G, H$ determine regular tetrahedra, respectively, and the five vertices $F, G, H, A, E$ define a regular triangular bipyramid. Finally, we complete the construction by adding the faces incident with either the north pole $A$ or the south pole $E$ (the faces are taken from the diagram of Figure 3, left). Since it can be easily checked with help from Figure 3 (left) that the graph's cycles (in the equatorial graph) around the north pole and south pole both are without self-intersections, it follows that the whole polyhedron constructed is free from self-intersections.

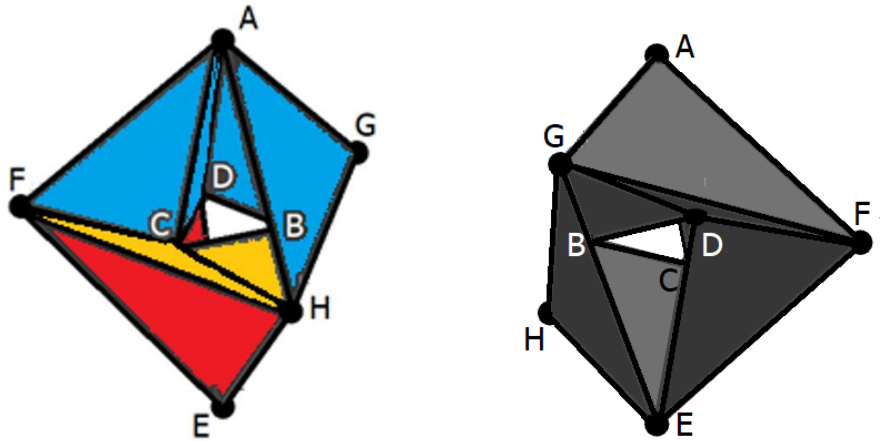

Figure 4: Toroidal polyhedral bipyramidal suspension at different angles, $k = 2.8$

A 3D model of the polyhedral suspension is constructed (with $k = 2.8$) with GeoGebra and is shown at two different typical angles in Figure 4.

Rotating the 3D model of the constructed suspension with GeoGebra, we accidentally came across with an unexpectedly elegant image shown in Figure 5.

Note that the geometric realization of the graph $O_4 = K_{2,2,2,2}$ described in this section only accommodates the triangulation in Figure 2, left (which is the same as Figure 1, left) as a geometric realization. For if we try to add the faces of the triangulation of Figure 2 (right) to the geometric graph, then faces $FDC$, $BDC$, and $BGH$ would entirely lie withing the face $FGH$ (check with Figures 2, 3), which intersections are impossible in a geometric realization of a triangulation. In Section 4 we construct another geometric realization of the graph $K_{2,2,2,2}$ which allows all the twelve triangulations with the vertex labeled graph $K_{2,2,2,2}$ as geometric realizations.

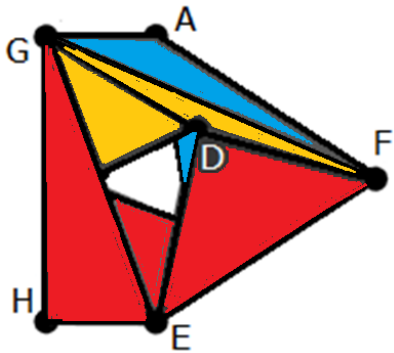

Figure 5: Designer variant for a logo

## 3   SCHLEGEL DIAGRAM: AN ENLIGHTENING EXAMPLE

In this section we provide a simple introductory example for a better understanding of the use of Schlegel diagrams in the forthcoming sections. Suppose we are given a usual 3D regular convex octahedron in $\mathbf{R}^3$ with the graph $O_3 = K_{2,2,2}$ and asked to geometrically realize its 1-skeleton (graph) in $\mathbf{R}^2$.

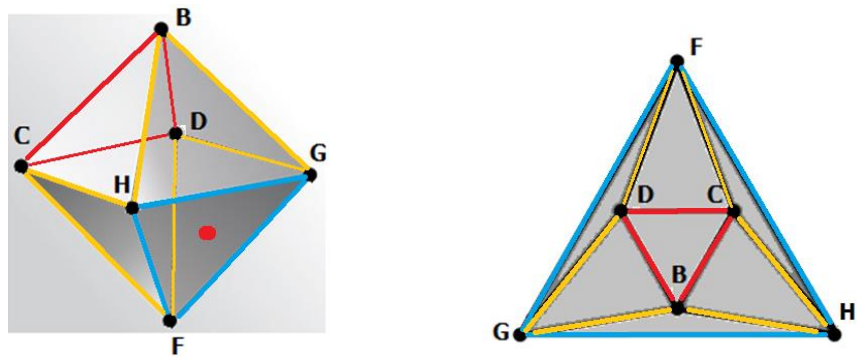

Figure 6: The octahedron and its Schlegel diagram

More generally, suppose we are given a convex $n$-polytope in $\mathbf{R}^n$ and asked to geometrically realize its $(n-2)$-skeleton in $\mathbf{R}^{n-1}$. A *Schlegel diagram* of a convex $n$-polytope in $\mathbf{R}^n$ is defined to be the central projection of the $n$-polytope from $\mathbf{R}^n$ into $\mathbf{R}^{n-1}$, onto one of its *facets* (that is, $(n-1)$-faces) through a point just outside the chosen facet (such a point is shown red in Figure 6 (left)). Due to the convexity of the $n$-polytope, it is not hard to see that the whole $(n-2)$-skeleton of the $n$-polytope is geometrically realized in the Schlegel diagram in $\mathbf{R}^{n-1}$.

Now, to solve the problem posed in the very beginning of this section, we take the Schlegel diagram of the regular convex octahedron as shown in Figure 6 (right), which provides a desired geometric realization of the whole graph of the octahedron in $\mathbf{R}^2$. Note that the inner ($BDC$) and outer ($GFH$) triangles in Figure 6 (right) are both regular and also note that the homothety with respect to the common circumcenter of the triangles, with some coefficient $-1/k$ ($k > 2$), maps the outer triangle onto the inner one.

# 4 NEW TOROIDAL POLYHEDRON

In this section we construct a polyhedron which is another geometric realization (in $\mathbf{R}^3$) of the abstract triangulation of the torus shown in Figures 1, 2 (left), which gives a new toroidal polyhedron in $\mathbf{R}^3$. Exploiting the 4D hyperoctahedron, our construction is a 4D generalization of the 3D method which uses the 3D octahedron ("8-cell") as outlined in Section 3. We use the 4D hyperoctahedron instead of the 3D octahedron. In fact, the hyperoctahedron contains, in its 2-skeleton, all the 12 triangulations of the torus with the vertex labeled graph $O_4 = K_{2,2,2,2}$ (this graph is the 1-skeleton of the hyperoctahedron); in particular it contains the pair of mutually complementary triangulations shown in Figure 2. We only need to project those 12 triangulations from $\mathbf{R}^4$ into one of the 3D facets (boundary tetrahedra) of the hyperoctahedron or in other words we need to realize the 12 triangulations geometrically in a Schlegel diagram of the hyperoctahedron.

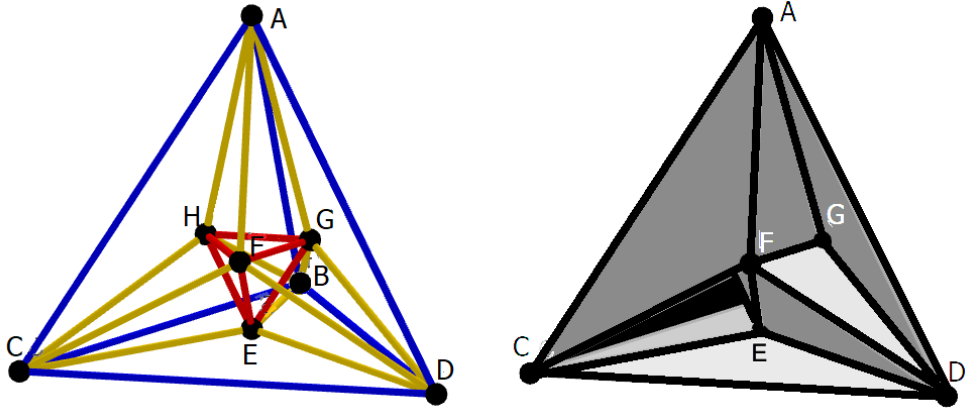

Figure 7: A geometric realization of the 1-skeleton (left)
and a view of the polyhedron with the interior of the face $ACD$ removed (right)

A Schlegel diagram is a projection of the hyperoctahedron (which is a convex 4-polytope) from $\mathbf{R}^4$ into $\mathbf{R}^3$ likewise we projected the usual octahedron (which is a convex 3-polytope) from $\mathbf{R}^3$ into $\mathbf{R}^2$ in Section 3. We extend the construction of Section 3 as follows (check with Figure 7). Fix on a regular tetrahedron $ABCD$ (blue) in $\mathbf{R}^3$ (check with Figure 7), whose circumcenter is at $(0,0,0)$. Also fix on its homothetic image $EFGH$ (red) with respect to $(0,0,0)$ with a negative coefficient $-1/k$ ($k > 0$) chosen so that $EFGH$ lies strictly inside $ABCD$. Then we connect the two tetrahedra by all possible edges with the exception of $AE$, $BF$, $CG$, and $DH$. The graph obtained this way is $O_4 = K_{2,2,2,2}$. Finally, we add the 16 faces of the triangulation as indicated in Figure 1 (left). In the following two paragraphs we assign coordinates to the 8 vertices in $\mathbf{R}^3$ and determine the range of $k$.

The coordinates of the vertices can be assigned as follows (this choice is not unique):

$$A = (0,0,3), \quad B = (\sqrt{8},0,-1), \quad C = (-\sqrt{2},\sqrt{6},-1), \quad D = (-\sqrt{2},-\sqrt{6},-1),$$
$$E = \left(0,0,-\frac{3}{k}\right), \quad F = \left(-\frac{\sqrt{8}}{k},0,\frac{1}{k}\right), \quad G = \left(\frac{\sqrt{2}}{k},-\frac{\sqrt{6}}{k},\frac{1}{k}\right), \quad H = \left(\frac{\sqrt{2}}{k},\frac{\sqrt{6}}{k},\frac{1}{k}\right) \quad (1)$$

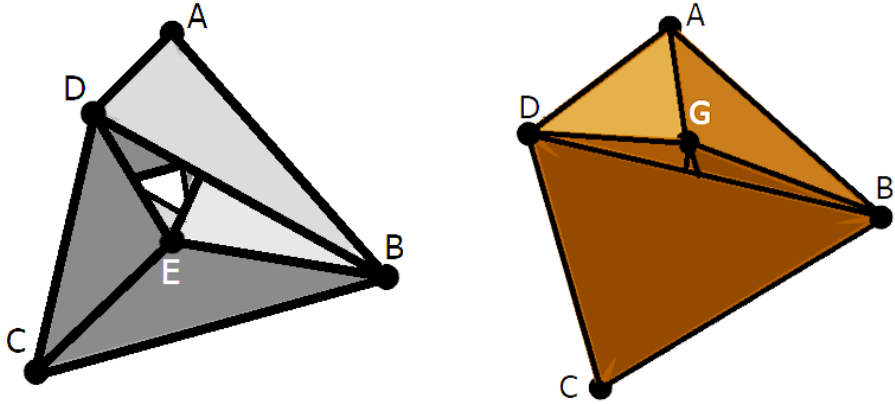

Figure 8: Left: A geometric realization of the abstract triangulation of Figures 1, 2 (left).
Right: A geometric realization of the abstract triangulation of Figure 2 (right), $k=4$

Now we determine $k$ for which the vertices of the inner tetrahedron (red) lie on the sphere inscribed into the outer tetrahedron (blue). The radius of the sphere circumscribed around the outer tetrahedron is 3, while the radius of the sphere circumscribed around the inner tetrahedron is equal to $3/k$. To find the value of $k$ so that the latter sphere turns out to be inscribed into the outer tetrahedron, note first that the side of the outer tetrahedron has length $2\sqrt{6}$, so the radius of the sphere inscribed into the outer tetrahedron is equal to 1. To have the latter sphere inscribed into the outer tetrahedron, its radius has to be equal to 1, so that $3/k=1$ and so $k=3$. Thus, $k$ has to be greater than 3 if we want the inner tetrahedron to lie entirely within the outer tetrahedron.

Thus, the vertices $A, B, C, D, E, F, G, H$ define a Schlegel diagram of the hyperoctahedron in $\mathbf{R}^3$ with coordinates of the eight vertices as given in Eqs. (1), and so the whole 2-skeleton of the hyperoctahedron is geometrically realized in $\mathbf{R}^3$. We come to the following.

**Theorem 1.** *All the twelve triangulations of the 2-torus with the vertex-labeled general octahedral graph $O_4 = K_{2,2,2,2}$ are realized geometrically with the same 1-skeleton in 3-space, with coordinates of the eight vertices as given in Eqs. (1).* ∎

We have implemented the construction, described in this section, with GeoGebra, taking $k=4$. In particular, Figure 8 presents 3D models of the new toroidal polyhedra which are geometric realizations of the abstract triangulations shown in Figure 2, respectively. We added the faces one by one, checking that each newly added face has no intersections with the faces already added, which made GeoGebra a kind of tool of experimental mathematics. Interestingly, Figure 8 presents two polyhedral tori (both without self-intersections) with the same 1-skeleton in 3-space, with the property that their intersection (as point-sets) is only their common 1-skeleton.

**Corollary 1.** *There exist six pairs of geometric toroidal polyhedra with the same 1-skeleton in 3-space, without a single common face in each pair.* ∎

## 5   GEOMETRIC REALIZATION OF THE TRIANGULATIONS OF THE 2D PROJECTIVE PLANE WITH THE COMPLETE GRAPH $K_6$

The coordinates of the six vertices in $\mathbf{R}^4$ for the polyhedron geometrically realizing the abstract triangulation of the 2D projective plane shown in Figure 1 (right) can be chosen as follows:

$$A = \left(0, 0, 0, \frac{4}{\sqrt{5}}\right), \qquad B = \left(1, 1, 1, -\frac{1}{\sqrt{5}}\right), \qquad C = \left(1, -1, -1, -\frac{1}{\sqrt{5}}\right),$$

$$D = \left(-1, 1, -1, -\frac{1}{\sqrt{5}}\right), \qquad E = \left(-1, -1, 1, -\frac{1}{\sqrt{5}}\right), \qquad O = (0, 0, 0, 0) \qquad (2)$$

Then the five vertices $A, B, C, D, E$ lie on the 3-sphere with radius $4/\sqrt{5}$, centered at the origin occupied by the sixth vertex $O = (0,0,0,0)$. The distance between any pair of the five vertices is $2\sqrt{2}$. Thus, the five vertices form a regular 4-simplex in $\mathbf{R}^4$, stellar subdivided by the sixth vertex $O = (0,0,0,0)$. Thus, the six vertices together form a Schlegel diagram of the 5-simplex in $\mathbf{R}^4$. Thus, the whole 2-skeleton of the 5-simplex is geometrically realized in $\mathbf{R}^4$ with coordinates of the six vertices as given in Eqs. (2). It follows that all the 12 triangulations with the vertex labeled graph $K_6$ are realized in $\mathbf{R}^4$ geometrically, with the same geometric graph, and we come to the following.

**Theorem 2.** *All the twelve triangulations of the 2D projective plane with the vertex labeled complete graph $K_6$ are realized geometrically with the same 1-skeleton in 4-space, with coordinates of the six vertices as given in Eqs. (2).* ∎

**Corollary 2.** *There exist six pairs of geometric 2D projective-planar polyhedra with the same 1-skeleton in 4-space, without a single common face in each pair.* ∎

# 6 GEOMETRIC REALIZATION OF THE TRIANGULATIONS OF THE MÖBIUS BAND WITH THE COMPLETE GRAPH $K_5$

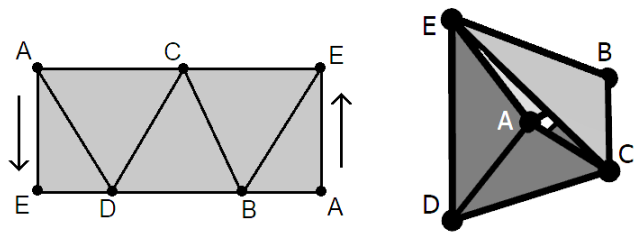

Figure 9: The initial abstract triangulation of the Möbius band with 5 vertices (left) and its geometric realization in 3-space (right)

In this section we obtain another corollary of Theorem 2. For that, we orthogonally project the geometric triangulation (as stated in Theorem 2) that realizes the abstract triangulation of the 2D projective plane (shown in Figure 1, right) into the $xyz$-hyperplane ($=\mathbf{R}^3$). Note that this projection contains a geometric realization (Figure 9, right) of an abstract triangulation of the Möbius band with five vertices, ten edges, and five faces (Figure 9, left). The five vertices are geometrically represented in $\mathbf{R}^3$ by the points with integer coordinates as follows:

$$A=(0,0,0,), \quad B=(1,1,1), \quad C=(1,-1,-1), \quad D=(-1,1,-1), \quad E=(-1,-1,1) \tag{3}$$

The abstract triangulation of the Möbius band shown in Figure 9 (left) is a subcomplex of the triangulation of the projective plane shown in Figure 1 (right). This triangulation appears to be the minimal triangulation of the Möbius band and is identified in [4] as one of the six so-called irreducible triangulations of the Möbius band (see [4] for details). As a corollary of Theorem 2 we have the following.

**Corollary 3.** *All the twelve triangulations of the Möbius band with the vertex labeled complete graph $K_5$ are realized geometrically with the same 1-skeleton in 3-space, with coordinates of the five vertices as given in Eqs. (3).* ∎

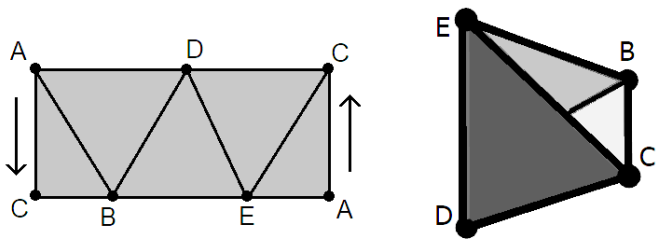

Figure 10: The complementary triangulation of the Möbius band with 5 vertices (left) and its geometric realization in 3-space (right)

In particular, Figure 10 (left) shows the abstract triangulation complementary to the triangulation in Figure 9 (left); they have the same vertex labeled graph but their face sets are disjoint.

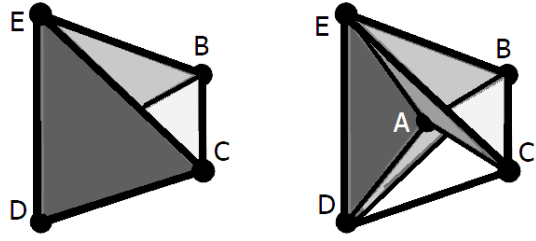

Figure 11: Geometric complementary triangulation of the Möbius band in 3-space (left) and the view of it with face *EDC* removed (right)

The left image of Figure 11 duplicates the geometric realization in the right image of Figure 10 while the right image of Figure 11 shows this geometric realization with the face *EDC* removed for unobstructed observation.

For a 4-simplex in $\mathbf{R}^4$, its Schlegel diagram is formed by a geometric 3-simplex (tetrahedron in $\mathbf{R}^3$) stellar subdivided by the 5th vertex chosen arbitrarily inside the 3-simplex.

**Corollary 4.** *The 2-skeleton of the 4-simplex is the face-disjoint union of two triangulated Möbius bands. Furthermore, both of them are realized geometrically in the Schlegel diagram of the 4-simplex in 3-space, which provides a pair of geometric polyhedral Möbius bands in 3-space with the same 1-skeleton, but without a single common face.* ∎

The geometric triangulation of the Möbius band, stated in Corollary 3, has edge lengths taking only two values, $2\sqrt{2}$ and $\sqrt{3}$, and has two equilateral and three isosceles triangles. One interesting open question is whether there exists a polyhedral geometric realization of the Möbius band in $\mathbf{R}^3$ with all edges of equal length. We conclude with another related open geometric question: If the projective plane, embedded in 4-space $xyzw$, is orthogonally projected into the hyperplane $xyz$, may the image of the projection be a Möbius band embedded in $\mathbf{R}^3$?

## 7  CONCLUSION

In the future we plan to extend the current research in five directions as follows:
(i) geometric symmetries [15],
(ii) volumes of polyhedra [16],
(iii) representations of the Cayley graph of the quaternion group $Q_8$ [15],
(iv) the relationship between 2D, 3D, and 4D,
(v) the inverse Newtonian potential problem for polyhedra.

Serge Lawrencenko
Russian State University of Tourism and Service, Institute of Service Technologies,
99 Glavnaya Street, Cherkizovo, Pushkino District, Moscow Region, 141221, Russia
lawrencenko@hotmail.com

Alex Lao
Kometa.Games, 6 Presnenskaya Naberezhnaya, Moscow, 123100, Russia
laoshanda@hotmail.com